\documentclass[12pt]{article}

\usepackage{amssymb}
\usepackage{amstext}
\usepackage{amscd}

\newtheorem{teo}{Th\'eor\`eme}[section]
\newtheorem{prop}[teo]{Proposition}
\newtheorem{lem}[teo]{Lemme}

\newtheorem{defini}[teo]{D\'efinition}

\newtheorem{rem}[teo]{Remarque}

\newcommand{\iii}{\sqrt{-1}}

\newcommand{\CC}{{\mathbb C}}
\newcommand{\RR}{{\mathbb R}}

\newcommand{\QQ}{{\mathbb Q}}
\newcommand{\NN}{{\mathbb N}}

\newcommand{\Hcal}{{\cal H}}
\newcommand{\GG}{{\mathbb G}}
\newcommand{\SSS}{{\mathbb S}}
\newcommand{\AAA}{{\mathbb A}}

\def\Fg{\mathfrak{g}}
\def\Fh{\mathfrak{h}}

\title{Equidistribution de sous-vari\'et\'es sp\'eciales.}
\author{ Laurent Clozel, Emmanuel Ullmo}

\begin{document}
\date{}
\maketitle


\section{Introduction}

Soit $S$ une vari\'et\'e de Shimura sur $\CC$.
On d\'efinit sur $S$ un ensemble de points sp\'eciaux
(les points \`a multiplication complexe) et un ensemble
de sous-vari\'et\'es sp\'eciales que l'on appelle sous-vari\'et\'es
de type de Hodge. Les d\'efinitions qui seront donn\'ees plus tard dans 
le texte sont pr\'esent\'ees de mani\`ere tr\`es agr\'eable dans le 
papier de Moonen \cite{Mo}.

Dans ce cadre Andr\'e et Oort font la conjecture suivante.
Soit $Y$ une sous-vari\'et\'e de $S$, il existe un ensemble
fini $\{S_{1},\ldots,S_{r}\}$ de sous-vari\'et\'es sp\'eciales
avec $S_{i}\subset Y$ pour tout $i$ tel que toute vari\'et\'e 
sp\'eciale $Z$ de $S$ contenue dans $Y$ est en fait contenue dans
un des $S_{i}$.   Le r\'esultat le plus
profond dans la  direction de cette conjecture a \'et\'e obtenu
par Edixhoven et Yafaev \cite{EdYa}.

On d\'efinit dans ce texte une classe assez large de sous-vari\'et\'es
sp\'eciales que nous appellerons fortement sp\'eciales par manque
d'une terminologie plus ad\'equate. D\'ecrivons 
les sous-vari\'et\'es fortement sp\'eciales:

Soit $S$  une vari\'et\'e de Shimura associ\'ee
\`a une donn\'ee de Shimura $(G,X)$ pour un groupe alg\'ebrique 
adjoint sur $\QQ$ et une $G(\RR)$-classe de conjuguaison $X$ de 
morphismes:
$$
h:\SSS \longrightarrow G_{\RR},
$$
o\`u $\SSS$ d\'esigne le tore de Deligne $Res_{\CC/\RR}\GG_{m}$.
Une sous-vari\'et\'e sp\'eciale de $S$ est associ\'ee \`a un 
$\QQ$-sous-groupe alg\'ebrique r\'eductif $H$. Les sous-vari\'et\'es
fortement sp\'eciales seront celles qui sont associ\'ees \`a un 
$\QQ$-sous-groupe alg\'ebrique semi-simple $H_{\QQ}е$ 
qui n'est contenu dans aucun $\QQ$-sous-groupe parabolique propre
de $G_{\QQ}$.е
 Le r\'esultat principal de 
ce texte est

\begin{teo} \label{teo1}
	 Soit $Y$ une sous-vari\'et\'e d'une vari\'et\'e de 
Shimura $S$.
Il existe un
   ensemble fini $\{S_{1},\ldots,S_{k}\}$ 
de sous-vari\'et\'es fortement sp\'eciales de dimension positive  
$S_{i}\subset 
Y$ tel que si $Z$ est une sous-vari\'et\'e fortement sp\'eciale de dimension
positive avec $Z\subset Y$ alors $Z\subset S_{i}$ pour un certain 
$i\in \{1,\dots,k\}$.
\end{teo}


Le th\'eor\`eme \ref{teo1} se d\'eduit d'un \'enonc\'e ergodique.
Toute sous-vari\'et\'e sp\'eciale $Z$ de $S$ est munie 
d'une mani\`ere canonique d'une mesure de probabilit\'e $\mu_{Z}$.

\begin{teo}\label{teo2}
Soit $S_{n}$ une suite de sous-vari\'et\'es fortement sp\'eciales.
Soit $\mu_{n}$ la mesure de probabilit\'e associ\'ee \`a $S_{n}$.
Il existe une sous-vari\'et\'e fortement sp\'eciale $Z$ et une
sous-suite $\mu_{n_{k}}$ qui converge faiblement vers $\mu_{Z}$.
De plus $Z$ contient $S_{n_{k}}$ pour tout $k$ assez grand.
\end{teo}
On obtient la preuve du th\'eor\`eme \ref{teo1} en consid\'erant
une suite de sous-vari\'et\'es fortement sp\'eciales maximales $S_{n}$
parmi les sous-vari\'et\'es fortement sp\'eciales contenues dans $Y$.
En passant \`a une sous-suite on peut supposer que $\mu_{n}$
converge faiblement vers $\mu_{Z}$. Comme le support de $\mu_{Z}$
est contenu dans $Y$, on en d\'eduit que $Z\subset Y$.
Par la maximalit\'e des $S_{n}$ et le fait que $S_{n}\subset Z$ pour
tout $n$ assez grand, on en d\'eduit que la suite $S_{n}$ est 
stationaire.

\bigskip

La preuve des r\'esutats principaux de ce texte repose sur des
r\'esultats ergodiques . L'outil principal de ce
texte est la conjecture de Raghunathan sur les flots unipotents
d\'emontr\'ee par Ratner \cite{Rat1} \cite{Rat2} et pr\'ecis\'ee par
Mozes et Shah \cite{MoSh}.  Dans la deuxi\`eme
partie de ce texte nous expliquons, dans le cadre 
arithm\'etique qui nous concerne, les  r\'esultats ergodiques
dont nous avons besoin. La troisi\`eme partie repose
essentiellement sur la th\'eorie des donn\'ees de Shimura $(G,X)$
d\'evelopp\'ee par Deligne \cite{De1} \cite{De2} interpr\'etant
les travaux de Shimura. On y montre les r\'esultats pr\'eliminaires
\`a la d\'emonstration des propri\'et\'es de stabilit\'e de l'ensemble
des sous-vari\'et\'es fortement sp\'eciales obtenues en d\'ebut de
quatri\`eme partie.  Les th\'eor\`emes principaux sont alors 
d\'emontr\'es \`a la fin de la quatri\`eme partie.  Nous donnons aussi
des exemples o\`u le th\'eor\`eme \ref{teo2} est mis en d\'efaut pour
des suites de vari\'et\'es sp\'eciales associ\'ees \`a des groupes $H_{n}$
qui ne sont pas semi-simples ou qui sont contenus dans un 
$\QQ$-parabolique propre.е

{\bf Remerciements.} Les auteurs remercient le rapporteur
pour d'utiles commentaires qui ont conduit \`a  une am\'elioration 
notable du r\'esultat principal de ce 
texte.

\section{Pr\'eliminaires sur les groupes}\label{section2}

{\bf Notations.}
Soit $H$ un groupe alg\'ebrique;  conform\'ement \`a l'usage 
on notera $H^0$ la composante connexe de $H$ pour la topologie 
de Zariski. $H^{ad}$,  $H^{der}$ et $H^{sc}$ 
d\'esignent respectivement le groupe adjoint, le groupe d\'eriv\'e
et le rev\^etement simplement connexe de $H^{der}$. 
On notera $R_{u}(H)$ le radical unipotent de $H$.
Si $H$ est un sous-groupe de $G$, on notera $N_{G}(H)$ le 
normalisateur dans $G$ de $H$ et $Cent_{G}(H)$ ou $Z_{G}(H)$ son 
centralisateur.  Si $H$ est semi-simple connexe et d\'efini sur un corps 
$k$,
$H$ est produit presque direct de ses $k$-sous-groupes connexes 
normaux minimaux $H_{1},\ldots,H_{r}ее$  (\cite{PR} prop2.4 p.62). 
Si $H$ est adjoint
ou simplement connexe ce produit est direct (\cite{PR} p.62).
Par abus de langage les $H_{i}$ seront appel\'es  facteurs $k$-simples
de $H$ dans la suite du texte.

Si $H_{1}$ est un facteur $\RR$-simple
d'un groupe semi-simple connexe $H$ \linebreak
sur $\RR$, on dit que $H_{1}$ est compact
ou non compact si $H_{1}(\RR)$ est compact ou non compact.
 On notera dans cette situation
$H(\RR)^+$ la composante connexe de $H(\RR)$ pour la topologie 
r\'eelle et 
$H(\RR)_+$
 la pr\'eimage de $H^{ad}(\RR)^+$ par l'application adjointe.
 Si  de plus $H$ est d\'efini sur $\QQ$, on note \linebreak
 $H(\QQ)^+ =H(\RR)^+\cap H(\QQ)$ et $H(\QQ)_+ =H(\RR)_+\cap H(\QQ)$.
Si $A$ est un sous-ensemble d'un espace topologique, on note
$\overline{A}$ son adh\'erence.

\bigskip 

Soient $G_{\QQ}$ un groupe alg\'ebrique connexe et semi-simple 
d\'efini sur $\QQ$ et
$G=G_{\QQ}(\RR)^+$. On suppose que  les groupes de
points r\'eels des facteurs 
$\QQ$-simples de $G_{\QQ}$ ne sont pas compacts. Soit $\Gamma$
un sous-groupe arithm\'etique de $G$ et $\Omega=\Gamma\backslash G$.
On note $P(\Omega)$ l'ensemble des mesures de Borel de probabilit\'e
sur $\Omega$.

Soit ${\cal H}$ l'ensemble des sous-groupes de Lie ferm\'es connexes
$H$ de $G$ tels que:

1) $H\cap \Gamma$ est un r\'eseau de $H$. En particulier 
$\Gamma\backslash \Gamma H$ est ferm\'e et on note $\mu_H\in P(\Omega)$
sa mesure $H$-invariante normalis\'ee.

2)  Le sous-groupe $L$ de $H$ engendr\'e par les sous-groupes 
unipotents \`a un param\`etre de $G$ contenus dans $H$ agit 
ergodiquement sur $\Gamma\backslash \Gamma H$  par rapport \`a
$\mu_H$.

Pour $H\in \Hcal$, on notera $L(H)$ (ou $L$ si il n'y a pas de 
confusion possible) le sous-groupe  de $H$ engendr\'e par les sous-groupes 
	 unipotents \`a un param\`etre de $G$ contenus dans $H$.

\begin{lem} \label{lem1}
	 Soient $H\in \Hcal$ et $L=L(H)$ le sous-groupe associ\'e.
	 
	 a) Soit $\overline{\Gamma\backslash \Gamma L}$
	 l'adh\'erence de $\Gamma\backslash 
	 \Gamma L$ dans $\Gamma\backslash G$.
	 Alors $\overline{\Gamma\backslash \Gamma L}=\Gamma\backslash 
	 \Gamma H$.

	 b)   Dans cette situation
	  $H$ est le plus petit sous-groupe de Lie ferm\'e de $G$ tel que 
$\overline{\Gamma\backslash \Gamma L}=\Gamma\backslash 
	 \Gamma H$.

	 c) Il existe un $\QQ$-sous-groupe alg\'ebrique $H_{\QQ}$
	 de $G_{\QQ}$ tel que $H=H_{\QQ}(\RR)^+$.
\end{lem}	 

{\it Preuve.} Notons tout d'abord que d'apr\`es les travaux de Ratner 
\cite{Rat1}   \cite{Rat2}, il existe un plus petit sous-groupe de Lie 
ferm\'e $H'$ de $G$ tel que $L\subset H'$ et 
$\overline{\Gamma\backslash \Gamma L}=\Gamma\backslash 
	 \Gamma H'$. D'apr\`es \cite{MoSh} (prop 2.1), $H'\in \Hcal$.
	 
	 Par ailleurs $L$ 
	 est un sous-groupe normal de $H$ et agit ergodiquement sur 
$\Gamma\backslash 
	 \Gamma H$.  Il existe donc une orbite sous $L$ qui est dense.
	 Il existe donc $h\in H$ tel que
	 
	 $$
	 \Gamma\backslash 
		  \Gamma H=\overline{\Gamma\backslash \Gamma h L}
		  =\overline{\Gamma\backslash \Gamma hLh^{-1}h}=
		  \overline{\Gamma\backslash \Gamma L}h.
	 $$

	 On en d\'eduit que $\overline{\Gamma\backslash \Gamma L}=\Gamma\backslash 
		  \Gamma H$; ce qui prouve (a).  
		  
		  Par minimalit\'e de $H'$, on a $H'\subset H$. 
		   D'apr\`es \cite{Sh}
		  (prop 3.2) si $H_{\QQ}$ d\'esigne le plus petit $\QQ$-sous-groupe
		  de $G_{\QQ}$ tel que $L\subset H_{\QQ}(\RR)$, on a 
		  $H_{\QQ}(\RR)^+=H=H'$. Ceci prouve donc (b) et (c).

		  \bigskip
		  
		  Si $E$ est un sous-ensemble de $G$, on d\'efinit le groupe de 
		  Mumford-Tate de $E$, not\'e $MT(E)$, comme le plus petit 
		  $\QQ$-sous-groupe alg\'ebrique $H_{\QQ}$ de $G_{\QQ}$ tel que
		  $E\subset H_{\QQ}(\RR)$. Si $H\in\Hcal$ et $L=L(H)$ 
		  alors $H=MT(L)(\RR)^+$. On retiendra le lemme suivant
		  d\^u \`a Shah.
	
		  \begin{lem}(Shah)\label{lemme2.2}
		  Soient $H\in \Hcal$ et $L=L(H)$. 
		  
		  a) Le radical $N$ de $L$ est unipotent et $L$ est un produit 
		  semi-direct
		  $$L=NS$$
		  pour un groupe semi-simple sans facteurs
		  compacts $S$.
		  
		  b) Le radical de $MT(L)$ est unipotent.
		     \end{lem}
		  {\it Preuve.} Le (a) est d\'emontr\'e dans \cite{Sh} Lemme 2.9.
		  Le (b) d\'ecoule de 
		  \cite{Sh} Prop. 3.2 et du fait que $\Gamma$ est un r\'eseau 
		  arithm\'etique (c.f. \cite{Sh}  Remarque 3.7).  
		  
		  \begin{lem}\label{lemme2.3}
		  Soit $H_{\QQ}$ un $\QQ$-sous-groupe alg\'ebrique connexe
		  semi-simple de $G_{\QQ}$.  Alors $H_{\QQ}(\RR)^+\in\Hcal$ 
		  si et seulement si 
		  pour tout facteur $\QQ$-simple  $H_{1\QQ}$ de  
		  $H_{\QQ}$, $H_{1\QQ}(\RR)$ n'est pas compact. 
		  \end{lem}
		  
		  {\it Preuve.}
		  Remarquons tout d'abord que par un r\'esultat de Cartan
		  (\cite{PR}, prop 7.6),
		  si $F$ est un $\RR$-groupe alg\'ebrique simple, simplement connexe et
		  non compact  alors $F(\RR)=F(\RR)^+$ 
		  est engendr\'e par ses sous-groupes
		  unipotents \`a un param\`etre. On en d\'eduit que si $F$
		  est un $\RR$-groupe alg\'ebrique simple non compact alors
		  $F(\RR)^+$ est engendr\'e par ses sous-groupes
		  unipotents \`a un param\`etre.

		Supposons que $H_{\QQ}$ est sans facteur $\QQ$-simple
		$\RR$-anisotrope.
		  Soit $L$ le sous-groupe de $H_{\QQ}(\RR)^+$ engendr\'e
		  par ses sous-groupes unipotents \`a un param\`etre.
		  Si $F$ est un facteur simple non compact de $H_{\QQ}(\RR)^+$,
		  alors par la discussion pr\'ec\'edente $F\subset L$.
		  On en d\'eduit que $MT(F)\subset MT(L)$. 
		  On en d\'eduit alors
		que $MT(L)$ contient les facteurs $\QQ$-simples de  $H_{\QQ}$ donc
		$MT(L)=H_{\QQ}$. 
		D'apr\`es les r\'esultats de Ratner
		(\cite{Rat3} thm 4 p. 162), il existe 
		$H'\in\Hcal$ minimal
		tel que $L\subset H'$ et $\Gamma\backslash \Gamma H'$ 
		soit ferm\'e dans $\Omega$. D'apr\`es le lemme \ref{lem1} on a
		 $H'=MT(L)(\RR)^+=H_{\QQ}(\RR)^+\in \Hcal$.
		  
		R\'eciproquement soit $H=H_{\QQ}(\RR)^+\in \Hcal$
		et $L=L(H)$.
		Si $H_{\QQ}е$ a un facteur $H_{1}$ $\QQ$-simple qui est $\RR$-anisotrope,
alors on a un morphisme surjectf
$$
\Gamma\cap H(\RR)^+\backslash H(\RR)^+\longrightarrow 
\Gamma_{1}\backslash 
H'_{1}(\RR)^+
$$
avec $H'_{1}$ isog\`ene \`a $H_{1}$, l'action de $L(H)$ \`a droite
\'etant triviale. L'image  $\Gamma_{1}$ de 
$\Gamma\cap H(\RR)^+$ est contenu dans un sous-groupe arithm\'etique (\cite{Bo2}, 
cor. 7.3) donc est finie. Ceci contredit l'ergodicit\'e de l'action de 
$L$.

	\subsection{Mesures alg\'ebriques.}	 
		
	Comme $G$ op\`ere \`a droite sur $\Omega$, on a une op\'eration
	induite de $G$ sur $P(\Omega)$ et pour $\mu\in P(\Omega)$, on note
	$\mu g$ son transform\'e par $g$.
	Soit $\mu\in P(\Omega)$, on note $\Lambda(\mu)$ son sous-groupe 
	d'invariance (donc ferm\'e dans $G$):
	$$
	\Lambda(\mu)=\{g\in G\ \ \vert \ \ \ \  \mu g=\mu\}
	$$
	et $Supp(\mu)$ son support. 
	On note $L(\mu)$ le sous-groupe de $G$ engendr\'e par les 
	sous-groupes \`a un param\`etre unipotents contenus dans 
	$\Lambda(\mu)$.
	On dit qu'une mesure $\mu\in P(\Omega)$ est alg\'ebrique si il existe
	$x\in \Omega$ tel que $Supp(\mu)=x\Lambda(\mu)$.
	
	On note $Q(\Omega)$ l'ensemble des $\mu\in P(\Omega)$ tels que 
	l'action de $L(\mu)$ sur $\Omega$ soit ergodique par rapport \`a $\mu$.
	D'apr\`es les r\'esultats de Ratner toute mesure dans $Q(\Omega)$
	est alg\'ebrique et d'apr\`es Mozes-Shah \cite{MoSh} pour tout 
	$\mu\in Q(\Omega)$, il existe un sous-groupe \`a un param\`etre
	unipotent $u(t)\in L(\mu)$ qui agit ergodiquement par rapport \`a 
	$\mu$. 	Le r\'esultat principal de \cite{MoSh} qui est \`a la base de
	ce texte est:
	
	\begin{teo}(Mozes-Shah) \label{Mozes-Shah}
	 Soit $\mu_i$ une suite de mesures dans $Q(\Omega)$ convergeant 
	vers $\mu\in P(\Omega)$.
	
	a) $Q(\Omega)$ est ferm\'e donc $\mu\in Q(\Omega)$. Soit
	$x\in supp(\mu)$.

	b) Soit $u_i(t)\subset L(\mu_i)$ un sous-groupe unipotent \`a un 
	param\`etre agissant ergodiquement par rapport \`a $\mu_i$. Soit 
	$g_i\in G$ une suite convergeant vers $e$ telle que $xg_i=x_i\in 
	supp(\mu_i)$ et telle que $\{xg_iu_i(t):\ \ t>0\}$ soit 
	\'equidistribu\'e par rapport \`a $\mu_i$
	(une telle suite existe \cite{MoSh}, p.156). Pour tout $i$ assez grand,
	on a 
	$$
	supp(\mu_i)\subset supp(\mu).g_i
	$$
	et
	$$
	g_i u_i(t) g_i^{-1}\in L(\mu).
	$$
	De plus le sous-groupe de $L(\mu)$ engendr\'e par les $g_i u_i(t) 
	g_i^{-1}$ pour $i$ assez grand agit ergodiquement par rapport \`a 
	$\mu$.
	\end{teo}

	En particulier soit $Q(\Omega,e)$, l'ensemble des mesures
	$\mu\in Q(\Omega)$ telles que $\Gamma .e\in supp(\mu)$.
	Les mesures de $Q(\Omega,e)$ sont les mesures $H$-invariantes
	normalis\'ees de support $\Gamma\backslash \Gamma H$ pour
	un $H\in \Hcal$. On utilisera aussi la proposition suivante
	essentiellement contenue dans Mozes-Shah \cite{MoSh}:

	\begin{prop}\label{prop2.5}
	L'ensemble $Q(\Omega,e)$ est compact
	pour la topologie faible. Si $\mu_{n}\in Q(\Omega,e)$
	est une suite  qui converge faiblement vers $\mu\in Q(\Omega,e)$, 
	alors pour
	tout $n$ assez grand $supp(\mu_{n}) \subset supp(\mu)$.
	\end{prop}

\section{Sous-vari\'et\'es sp\'eciales des vari\'et\'es de Shimura.}

\subsection{Pr\'eliminaires}\label{section3.1}

Soit ${\SSS}=Res_{\CC/\RR} \GG_{m\CC}$ le tore de Deligne, une 
donn\'ee de Shimura est un couple $(G_{\QQ},X)$ o\`u $G_{\QQ}$ est un groupe 
r\'eductif sur $\QQ$ et $X\in Hom(\SSS,G_{\RR})$ est une classe de 
$G(\RR)$-conjugaison v\'erifiant les ``conditions de Deligne''
\cite{De1}, \cite{De2}:

a) Pour tout $\alpha\in X$ la repr\'esentation adjointe $Lie(G_{\RR}$)
est de type $$\{(-1,1),(0,0),(1,-1)\};$$
en particulier
$\alpha(\GG_{m\RR})\subset Z(G_{\RR})$.

b) L'involution $int(\alpha(\iii))$ est une involution de Cartan du 
groupe adjoint $G_{\RR}^{ad}$.

c) Le groupe $G_{\QQ}^{ad}$ n'a pas de $\QQ$-facteur $\RR$-anisotrope.е

\medskip

On suppose dans la suite de cette section que $G_{\QQ}$ est adjoint. Pour tout 
$\alpha\in X$, le groupe de Mumford-Tate $MT(\alpha)$
est d\'efini comme le plus petit $\QQ$-sous-groupe de 
$G_{\QQ}$ tel que l'on ait une factorisation de $\alpha$
via $MT(\alpha)_{\RR}$. (Noter que ce groupe est donc connexe).
Quand $T=MT(\alpha)$ est un tore,
on dit que $\alpha$ est sp\'ecial; comme $T({\RR})$
est contenu dans le centralisateur de $\alpha(\iii)$ qui est 
compact, on en d\'eduit que $T({\RR})$ est compact.

\begin{defini}
 Une sous-donn\'ee de Shimura $(H_{\QQ}е,X_{H})$ de $(G_{\QQ}е,X)$ 
 est une donn\'ee
 de Shimura telle que $H_{\QQ}е$ est un $\QQ$-sous-groupe alg\'ebrique 
 de $G_{\QQ}е$
 et $X_{H}$ la $H(\RR)$-classe de conjugaison d'un morphisme
 $\alpha: \SSS \rightarrow G_{\RR}$,  $\alpha\in X$
 se factorisant par $H_{\RR}$.
	 \end{defini}	 

	 \begin{prop}\label{Cartan4} 
Soit $H_{\QQ}$ un $\QQ$-sous-groupe alg\'ebrique 
de $G_{\QQ}$ semi-simple connexe et
sans $\QQ$-facteur $\RR$-anisotrope. On suppose qu'il existe 
$\alpha:\SSS\rightarrow G_{\RR}$,  $\alpha\in X$ se factorisant par
$H_{\RR}$. Soit $X_{H}$ la $H(\RR)$-classe de conjuguaison de $\alpha$.е
 Alors $(H_{\QQ},X_{H})$ est une 
sous-donn\'ee de Shimura de $(G_{\QQ},X)$.
\end{prop}

Nous v\'erifions les conditions (a), (b) et (c) des donn\'ees de 
Shimura. Soit $\Fh=\mbox{Lie } H(\RR)$, $\Fg= \mbox{Lie } G(\RR)$,
$C=\alpha(\sqrt{-1})\in H(\RR)$; alors $C^2$ est central dans $H(\RR)$.
 La condition (a) d\'ecoule du fait que $\Fh$  est un sous-espace de $\Fg$ invariant
par $\SSS$. 

Pour (b), $H$ \'etant semi-simple, il nous suffit de v\'erifier que
$\mbox{int}(C)$ est une involution de Cartan de $H_{\RR}$.
D'apr\`es (\cite{De2}, 1.1.15), il suffit d'exhiber une 
repr\'esentation r\'eelle  $V$ de $H(\RR)$, fid\`ele et $C$-polarisable
au sens suivant: il existe une forme bilin\'eaire $B$ sur $V$, 
invariante,  telle que  $B(X,CY)$ soit sym\'etrique et d\'efinie 
positive. On prend $V=\Fg$ pour la repr\'esentation adjointe et $B$
\'egale \`a la forme de Killing.  Enfin (c) est vrai par hypoth\`ese.

\subsection{Sous-vari\'et\'es de type de Hodge}
On note $\AAA$ l'anneau des ad\`eles de $\QQ$ et $\AAA_f$ l'anneau 
des ad\`eles finis. Soit $(G,X)$ une donn\'ee de Shimura ($G$ 
n'\'etant pas n\'ecessairement  adjoint)
et $K$ un 
sous-groupe compact ouvert de $G(\AAA_f)$, on note
$$
Sh_K(G,X)(\CC)= G(\QQ)\backslash X\times G(\AAA_f)/K
$$
et $[x,gK]$ l'image de $(x,gK)\in X\times G(\AAA_f)$ dans
$Sh_K(G,X)(\CC)$.

Soit $X^+$ une composante connexe de $X$; $X^+$ est une 
$G^{ad}(\RR)^+$-classe de conjugaison d'un morphisme
$h^{ad} :\ \ \SSS\rightarrow G_{\RR}^{ad}$ et $X^+$ est un domaine
sym\'etrique hermitien.  Soit $K_{\infty}$ le fixateur de 
$h^{ad}(\iii)$ dans $G^{ad}(\RR)^+$. Soit $K_{\infty,+}$ la 
pr\'eimage de $K_{\infty}$ par l'application adjointe, on a alors un
isomorphisme
\begin{equation}\label{Shimura1}
X^{+}\simeq G(\RR)_{+}/K_{\infty,+}\simeq G^{ad}(\RR)^+/K_{\infty}.
\end{equation}
et 
\begin{equation}
Sh_K(G,X)(\CC)= G(\QQ)_{+}\backslash X^{+}\times G(\AAA_f)/K.
\end{equation}
On note encore $[x,gK]$ l'image de $(x,gK)\in X^{+}\times G(\AAA_f)$ dans
$Sh_K(G,X)(\CC)$.

Nous aurons besoin de la d\'efinition des op\'erateurs de Hecke 
dans ce cadre (voir par exemple \cite{Mo2} 1.6.1).
\begin{defini}
	 Soient $g\in G(\AAA_{f})$ et  $K_{g}=K\cap g K g^{-1}$. La 
	 correspondence de Hecke $T_{g}$ sur $Sh_K(G,X)(\CC)$ est d\'efinie
	 par le diagramme
	 $$
	 Sh_K(G,X)(\CC)\leftarrow^{\pi_{1}е}\ \  Sh_{K_{g}е}(G,X)(\CC) \ \ \ 
	 ^{\pi_{2}е}\rightarrow 
	 Sh_K(G,X)(\CC).
	 $$
	 o\`u $\pi_{1}$ est donn\'e par l'inclusion $K_{g}\subset K$ et
	 $\pi_{2}е$ est l'application
	 $$
	 [x,\theta]\rightarrow [x,\theta g].
	 $$
	Soit $Z$ une sous-vari\'et\'e de $Sh_K(G,X)(\CC)$, on note 
	$T_{g}.Z$ le cycle $\pi_{2*}\pi_{1}^*Z$ de $Sh_K(G,X)(\CC)$.
	On dit que $T_{g}. Z$ est le translat\'e de $Z$ par l'op\'erateur
	de Hecke $T_{g}$.
\end{defini}

Soit $R_{G,K}$ un syst\`eme de repr\'esentants de 
$G(\QQ)_{+}\backslash G(\AAA_f)/K$, alors $R_{G,K}$ est fini et
\begin{equation}\label{Shimura2}
Sh_K(G,X)=\cup_{g\in R_{G,K}} \Gamma_g\backslash X^+
\end{equation}
o\`u 
$$
\Gamma_g=G(\QQ)_+\cap gKg^{-1}.
$$ 
Si $\Gamma'_{g}$ d\'esigne l'image par l'application adjointe de 
$\Gamma_{g}$ on a un isomorphisme
$$
\Gamma'_g\backslash X^+=\Gamma_g\backslash X^+
$$
o\`u les groupes $\Gamma_{g}$ et $\Gamma'_{g}$ agissent de mani\`ere
naturelle via les isomorphismes de l'\'equation (\ref{Shimura1}).

On suppose dans la suite de cette section que $G=G^{ad}$ est un groupe 
adjoint donc que  $X^{+}$ est une $G(\RR)^{+}$ classe de conjugaison 
de morphismes de 
$$\SSS \rightarrow G_{\RR}$$
et  $\Gamma_{g}=G(\QQ)^{+}\cap gKg^{-1}$.

Soit $(H,X_{H})$ une sous-donn\'ee de Shimura. Si $K_{H}=K\cap 
H(\AAA_{f})$,
on dispose d'un morphisme induit de vari\'et\'es de Shimura
$$
\psi:\ Sh_{K_{H}}(H,X_{H})(\CC)\longrightarrow Sh_K(G,X)(\CC).
$$
On choisit alors un syst\`eme de repr\'esentant $R_{H,K}$ de 
$$
H(\QQ)_{+}\backslash H(\AAA_{f})/K_{H},
$$
on a donc 
$$
Sh_{K_{H}}(H,X_{H})(\CC)=
\cup_{\lambda\in R_{H,K}} \Delta_{\lambda}\backslash X_{H}^{+}
$$
avec $\Delta_{\lambda}=H(\QQ)_{+}\cap \lambda K_{H}\lambda^{-1}.$

\begin{defini} Avec les notations pr\'ec\'edentes,
une sous-vari\'et\'e de la forme $\psi(\Delta_{\lambda}\backslash 
X_{H}^{+})$ est appel\'ee sous-vari\'et\'e de 
type Shimura de $Sh_K(G,X)(\CC)$. Une composante
irr\'eductible d'un translat\'e par un op\'erateur 
de Hecke d'une sous-vari\'et\'e de type Shimura de $Sh_K(G,X)(\CC)$
est appel\'ee sous-vari\'et\'e de type de Hodge. 
\end{defini}

Le but de cette partie est de d\'ecrire les sous-vari\'et\'es de type
de Hodge dans le langage des espaces localement sym\'etriques
hermitiens. Le lemme suivant qui montre la faible diff\'erence entre
les notions de sous-vari\'et\'e de type Shimura et sous-vari\'et\'e 
de type de Hodge nous permettra de nous ramener toujours dans la suite
\`a des sous-vari\'et\'es de type Shimura.

\begin{lem}\label{lemme3.7}
Soit $M$ une sous-vari\'et\'e de type de Hodge de $Sh_K(G,X)(\CC)$.
Il existe $\beta\in R_{G,K}$ et une sous-vari\'et\'e de type Shimura
$M_{1}$ tels que $M$ est une composante irr\'eductible de 
$T_{\beta}.M_{1}$.
\end{lem}

{\it Preuve.} Il existe une sous-donn\'ee de Shimura $(H,X_{H})	$ et
$\lambda\in G(\AAA_{f})$ tels que $M$ est
l'image de $X_{H}\times \lambda K$ dans $Sh_K(G,X)(\CC)$.
On peut \'ecrire $\lambda=\gamma \beta k$ avec $\gamma\in G(\QQ)_{+}$,
$\beta\in R_{G,K}$ et $k\in K$. 
Soient $H_{\gamma}=\gamma^{-1} H\gamma$ et $X_{\gamma}$ 
la $H_{\gamma}(\RR)$-classe de conjugaison de $\gamma^{-1} .x_{0}$
pour un $x_{0}\in X_{H}$ , $(H_{\gamma},X_{\gamma})$ est
une sous-donn\'ee de Shimura et $M$ est aussi l'image 
de   $X_{\gamma}\times \beta K$  dans
$Sh_K(G,X)(\CC)$.  On en d\'eduit que $M$ est une composante 
irr\'eductible de  
$$T_{\beta}.
Sh_{K\cap H_{\gamma}(\AAA_{f})}(H_{\gamma},X_{\gamma})(\CC).$$

\begin{lem}\label{lemme2.6}
Pour $\lambda\in R_{H,K}$, il existe un unique $\beta\in R_{G,K}$
tel que $\lambda=\gamma \beta k$ avec $\gamma\in G(\QQ)^{+}$
et $k\in K$. On a alors
$$
\psi (\Delta_{\lambda}\backslash X_{H}^{+})\subset 
\Gamma_{\beta}\backslash X^+
$$
\end{lem}	 

{\it Preuve.} On a pour tout $x\in X^{+}_{H}$
$$
\psi([x,\lambda K_{H}])=[x,\lambda K]=
[x,\gamma \beta K]=[\gamma^{-1}x,\beta K].
$$
Ceci termine la preuve quand on a remarqu\'e que les \'el\'ements
de $\Delta_{\lambda}\backslash X_{H}^{+}$ sont ceux de la forme
$[y,\lambda K_{H}]$ ($y\in X_{H}^{+}$) et ceux de
$\Gamma_{\beta}\backslash X^+$ sont ceux
de la forme $[y,\beta K]$ 
avec $y\in X^+$.

\medskip

Fixons $x_{0}\in X_{H}^+$ de sorte que 
$$X_{H}^+=H(\RR)_+.x_{0}\subset X^{+}= G(\RR)^{+}.x_{0}.$$
Soient $x_{1}=\gamma^{-1}.x_{0}\in X$ et 
$H_{\gamma}=\gamma^{-1} H \gamma$. On a 
$H_{\gamma}(\RR)=\gamma^{-1} H(\RR) \gamma$ et on note
$$
X_{H_{\gamma}}=H_{\gamma}(\RR).x_{1} 
$$ 
la $H_{\gamma}(\RR)$-classe de conjugaison
de $x_{1}$ alors
$$
X_{H_{\gamma}}^{+}=H_{\gamma}(\RR)_{+}.x_{1}
$$
est une composante connexe de $X_{H_{\gamma}}$.

On note $\psi_{\lambda}$ l'inclusion naturelle
$$
\psi_{\lambda} : X_{H_{\gamma}}^+\longrightarrow X^{+}.
$$
\begin{lem}\label{lemme2.8}
a) L'application $\psi_{\lambda}$ induit par passage au quotient une 
application (encore not\'ee $\psi_{\lambda}$)
$$
\psi_{\lambda}: \     \gamma^{-1} \Delta_{\lambda}\gamma \backslash
X_{H_{\gamma}}^+\longrightarrow \Gamma_{\beta}\backslash X^+,
$$
et 
$$
\psi_{\lambda}(    \gamma^{-1} \Delta_{\lambda}\gamma \backslash
X_{H_{\gamma}}^+)=\psi(   \Delta_{\lambda}\backslash X_{H}^{+}  ).
$$
b) On a $\gamma^{-1} \Delta_{\lambda}\gamma \subset \Gamma_{\beta}$
et 
$$
\gamma^{-1} \Delta_{\lambda}\gamma =
H_{\gamma}(\QQ)_{+}\cap \Gamma_{\beta}=
H_{\gamma}(\RR)_{+}\cap \Gamma_{\beta}.
$$
\end{lem}
{\it Preuve.} Comme $\gamma^{-1}\lambda=\beta k$, on a 
$$
\gamma^{-1} \Delta_{\lambda} \gamma=H_{\gamma}(\QQ)_{+}
\cap \beta k K_{H}\beta^{-1}\subset \Gamma_{\beta}.
$$
Ceci prouve \`a la fois la premi\`ere partie du (a) et du (b).
Par ailleurs d'apr\`es la preuve du lemme \ref{lemme2.6}
$$
\psi (\Delta_{\lambda}\backslash X_{H}^{+})=\{[\gamma^{-1}h.x_{0}, 
\beta K],  \ \  \mbox{  } h\in H(\RR)_{+}\}
$$
d'o\`u 
$$
\psi (\Delta_{\lambda}\backslash X_{H}^{+})=
\{[h.x_{1}, 
\beta K], \mbox{ } h'\in H_{\gamma}(\RR)_{+}\}=
\psi_{\lambda}(    \gamma^{-1} \Delta_{\lambda}\gamma \backslash
X_{H_{\gamma}}^+).
$$
Comme $\Gamma_{\beta}\subset G(\QQ)$, on a 
$$
H_{\gamma}(\QQ)_{+}\cap \Gamma_{\beta}=
H_{\gamma}(\RR)_{+}\cap \Gamma_{\beta}.
$$
On a 
\begin{equation}\label{ZUT}
\gamma^{-1} \Delta_{\lambda}\gamma =H_{\gamma}(\QQ)_{+}\cap
\beta k K_{H} k^{-1} \beta^{-1}
\end{equation}
et
$$
H_{\gamma}(\QQ)_{+}\cap \Gamma_{\beta}=H_{\gamma}(\QQ)_{+}\cap 
\beta K \beta^{-1}.
$$
Un \'el\'ement $\theta\in      H_{\gamma}(\QQ)_{+}\cap \Gamma_{\beta} $
peut donc s'\'ecrire
$$
\theta= \beta k_{1}\beta^{-1}=\gamma^{-1} h \gamma 
$$
avec $k_{1}\in K$ et $h\in H(\QQ)_{+}$.  On a donc
$$
k^{-1} k_{1} k= \lambda^{-1} h \lambda\in H(\AAA_{f})\cap K=K_{H}
$$
et $\theta\in H_{\gamma}(\QQ)_{+}\cap
\beta k K_{H} k^{-1} \beta^{-1}$. Ceci termine la d\'emonstration du 
lemme au vu de l'\'equation (\ref{ZUT}).

R\'esumons l'information qui nous sera utile dans la suite sous la 
forme:

\begin{prop}\label{prop2.9}
	 On suppose toujours que $G_{\QQ}$ est adjoint et on fixe 
	 $\beta\in R_{G,K}$. On pose $\Gamma=\Gamma_{\beta}$ de 
	 sorte que $S_{0}=\Gamma\backslash X^+$ est une composante 
	 irr\'eductible de $Sh_{K}(G,X)(\CC)$. 
 
	Soit $M$  une sous-vari\'et\'e de type Shimura de
 $S_{0}$, alors $M$ est l'image dans $S_{0}$ d'une vari\'et\'e
$M'=\Delta_{H}\backslash X_{H}^{+}$ o\`u $H_{\QQ}$ est un 
 un $\QQ$-sous-groupe de $G_{\QQ}$
tel que :

1) $\Delta_{H}= H(\RR)_{+}\cap \Gamma$ est un r\'eseau
arithm\'etique de $H(\RR)_{+}$.

2) Il existe un sous-groupe 
compact maximal $K_{\infty}$ de $G(\RR)^{+}$ tel que \linebreak
$K_{\infty}\cap H(\RR)_{+}$ est un compact maximal de 
$H(\RR)_{+}$ et $$X_{H}^+ \simeq H(\RR)_{+}/ K_{\infty}\cap H(\RR)_{+}.$$
\end{prop}
On a aussi une r\'eciproque utile \`a cette proposition:
\begin{prop}\label{prop3.11}
Soit $H_{\QQ}$ un $\QQ$-sous groupe r\'eductif v\'erifiant
les $2$ propri\'et\'es de la proposition \ref{prop2.9} avec
$K_{\infty}=Cent(\alpha(\iii))$ pour un $\alpha\in X$ tel que 
$MT(\alpha)\subset H_{\QQ}е$ soit un tore, alors l'image $M$ de 
$\Delta_{H}\backslash X_{H}^{+}$ dans $S_{0}$ est une sous-vari\'et\'e
de type de Hodge.
\end{prop}
{\it Preuve.} La sous-vari\'et\'e $M$ est totalement g\'eod\'esique 
dans $S_{0}$ et contient un point sp\'ecial,  par les r\'esultats de 
Moonen (\cite{Mo} thm 4.3),
c'est une 
sous-vari\'et\'e de type de Hodge.

\section{Preuve des Th\'eor\`emes }

\subsection{Sous-vari\'et\'es fortement sp\'eciales}

Soient $(G,X)$ une donn\'ee de Shimura avec $G$ adjoint, $K$ 
un sous groupe compact ouvert de $G(\AAA_{f})$ et $S=Sh_{K}(G,X)(\CC)$.
Une sous-vari\'et\'e fortement sp\'eciale est une composante 
irr\'eductible d'un translat\'e par un op\'erateur de Hecke d'une 
sous-vari\'et\'e de Shimura $Sh_{K\cap H'(\AAA_{f})}(H',X_{H'})(\CC)$
o\`u 

a) $H'$  est semi-simple.

b) $H'$ n'est pas contenu dans un $\QQ$-sous-groupe
parabolique propre de $G_{\QQ}$.

D'apr\`es (\cite{EMS}, lemme 5.1),  
la condition (b) est \'equivalente 
aux conditions (b') ou  (b") suivantes:

b') Tout $\QQ$-sous-groupe alg\'ebrique de $G_{\QQ}$ contenant
$H_{\QQ}$ est r\'eductif

b") Le centralisateur $Z_{G}(H)$ de $H_{\QQ}$ dans $G_{\QQ}$
est $\QQ$-anisotrope.

\medskipее

Plus g\'en\'eralement si on ne suppose plus $G$ adjoint,
soient $(G,X)$ une donn\'ee de Shimura et $K\subset G(\AAA_{f})$ un
sous-groupe compact ouvert. 
Soit $(G^{ad},X^{ad})$ la donn\'ee de Shimura adjointe (\cite{Mo2} 
1.6.7). Soit $K^{ad}$ un sous-groupe compact ouvert contenant l'image
par l'application adjointe de $K$; le morphisme
induit de $S=Sh_{K}(G,X)(\CC)$ vers $S^{ad}=Sh_{K^{ad}}(G^{ad},X^{ad})(\CC)$ 
est fini et d'apr\`es \cite{EdYa} (proposition 2.2)
une sous-vari\'et\'e  $Z$ de $S$
est de type de Hodge si et seulement si son image $Z^{ad}$ dans
$S^{ad}$ l'est. On dit alors que $Z$ est fortement sp\'ecial si 
$Z^{ad}$ l'est.

On suppose de nouveau $G$ adjoint.
On fixe encore $\beta\in R_{G,K}$. On pose $\Gamma=\Gamma_{\beta}$
et $S_{0}=\Gamma\backslash X^+$. On rappelle que l'on a 
pu d\'efinir avec les m\^emes notations
dans la section 
\ref{section2} un ensemble de sous-groupes
de Lie connexes de $G(\RR)^{+}$
not\'e $\Hcal$.  Si $\Omega=\Gamma \backslash G(\RR)^+$, on
a aussi d\'efini des ensembles de mesures de probabilit\'es 
$$Q(\Omega,e) \subset Q(\Omega)\subset P(\Omega).$$

D'apr\`es la proposition 
\ref{prop2.9}, une sous-vari\'et\'e fortement sp\'eciale
de $S_{0}$ associ\'ee \`a une sous-donn\'ee de Shimura $(H',X_{H'})$
de $(G,X)$
est (\`a translation par un op\'erateur de Hecke pr\`es)
l'image d'une vari\'et\'e de la forme 
$\Delta_{H}\backslash X_{H}^+$ v\'erifiant les conditions de la 
proposition \ref{prop2.9} pour un groupe $H_{\QQ}$ qui est un 
conjugu\'e de $H'_{\QQ}$ par un \'el\'ement de $G(\QQ)$.
On en d\'eduit que $H$ v\'erifie les m\^emes 
propri\'et\'es a) et b) que $H'$. 
Par abus de notation, on \'ecrira souvent
$\Delta_{H}\backslash X_{H}^+$ pour son image dans $S_{0}$.

\begin{rem}\label{remarque}
{\rm Le sous-groupe $H_{\QQ}$, associ\'e \`a 
une sous-vari\'et\'e sp\'eciale $M=  \Delta_{H}\backslash X_{H}^+  $ 
n'est bien d\'efini qu'a
conjugaison pr\`es par un $\lambda\in \Gamma$. Si $X_{H}^+=H(\RR)_{+}.x_{0}$
pour un $x_{0}\in X_{H}$, on note $H_{\lambda}=\lambda H_{\QQ} \lambda^{-1}$,
$x_{1}=\lambda.x_{0}$, $X_{H_{\lambda}}^+$ la 
$H_{\lambda}(\RR)_{+}$-classe de conjugaison de $x_{1}$ et 
$\Delta_{H_{\lambda}}=\Gamma \cap H_{\lambda}(\RR)_{+}$.
Alors $H_{\lambda}$ a les m\^emes propri\'et\'es que $H$ et 
$M$ est aussi l'image de $\Delta_{H_{\lambda}}\backslash 
X_{H_{\lambda}}^+$.}
\end{rem}	 
Par ailleurs d'apr\`es le lemme
\ref{lemme2.3}  et la condition c) de Deligne pour les vari\'et\'es de 
Shimura (rappell\'e au d\'ebut de la section \ref{section3.1}), on a
$H(\RR)^{+}\in \Hcal$.

Soit $M=\Delta_{H}\backslash X_{H}^+$ 
une sous-vari\'et\'e fortement sp\'eciale de type Shimura. 
Soit $\alpha\in X^+$ se factorisant par $H_{\RR}$.е
S'il existe un sous-groupe de Lie connexe $F\in \Hcal$ tel que 
$H_{\QQ}(\RR)^+\subset F$, on sait (d'apr\`es le lemme \ref{lem1})
qu'il existe un $\QQ$-sous-groupe 
$H'_{\QQ}$ tel que $F=H'(\RR)^+$.  D'apr\`es la propri\'et\'e (b)
des sous-vari\'et\'es fortement sp\'eciales, $H'_{\QQ}$ est
r\'eductif (sinon $H'_{\QQ}$ donc $H_{\QQ}$ serait contenu
dans un $\QQ$-parabolique propre. 
Soit $X_{H'}^+$ la 
$H'(\RR)_{+}$-classe de conjugaison de $\alpha$ et 
$\Delta_{H'}=\Gamma\cap H'(\RR)_{+}$.

\begin{lem}\label{lemme4.1}
La  sous-vari\'et\'e $M'=	\Delta_{H'}\backslash  X_{H'}^+$ est
fortement sp\'eciale et 
$$M\subset M'.$$
\end{lem}

{\it Preuve. } On a vu que $H'_{\QQ}$ est r\'eductif. 
Comme $H'_{\QQ}(\RR)^+\in \Hcal$,  on sait d'apr\`es les lemmes
\ref{lemme2.2} et \ref{lemme2.3} que  $H'_{\QQ}$ est semi-simple
sans $\QQ$-facteur simple compact.е
D'apr\`es la proposition \ref{Cartan4}, $(H',X_{H'})$ 
est une sous-donn\'ee de Shimura, on en d\'eduit par la proposition
\ref{prop3.11} que $M'$ est une sous-vari\'et\'e de type de Hodge.
 L'assertion $M\subset M'$ est claire. 

On rappelle que pour tout $F\in \Hcal$ on dispose d'une mesure
de probabilit\'e associ\'ee $\mu_{F}$ sur 
$Z_{F}=\Gamma\backslash \Gamma F
\simeq \Gamma\cap F\backslash F$. En particulier si 
$M=\Delta_{H}\backslash X_{H}^+$ est une sous-vari\'et\'e fortement 
sp\'eciale de type Shimura, comme $H_{\QQ}(\RR)^+\in \Hcal$,
on dispose d'une mesure de probabilit\'e
$\mu_{H}=\mu_{H(\RR)^+}$ sur 
$$
Z_{H}=Z_{H(\RR)^+}=H(\RR)^+\cap \Gamma \backslash H(\RR)^+.
$$  

Par ailleurs 
soit $\alpha\in X^+$ tel que $T=MT(\alpha)\subset H_{\QQ}$ est un tore.
Alors $K_{\infty}=Cent(\alpha(\iii))$   est un compact maximal de 
$G(\RR)^+$,  
$$K_{\infty}\cap H(\RR)_{+}$$
est un compact maximal
de $H(\RR)_{+}$ et 
$$X_{H}^+=H(\RR)_{+}/K_{\infty}\cap H(\RR)_{+}.$$

Alors $K_{\infty}\cap H(\RR)^{+}$ est un compact maximal
de $H(\RR)^{+}$ et comme $K_{\infty}$ agit transitivement
sur $H(\RR)_{+}/H(\RR)^{+}$ on dispose d'un isomorphisme
$$
H(\RR)^+/K_{\infty}\cap H(\RR)^{+} \longrightarrow X_{H}^+.
$$
On en d\'eduit en composant les applications naturelles 
$$
 Z_{H}= H(\RR)^+\cap \Gamma \backslash H(\RR)^+\longrightarrow 
 H(\RR)^+\cap \Gamma \backslash  X_{H}^+ 
 $$
 et
 $$
 H(\RR)^+\cap \Gamma \backslash  X_{H}^+ \longrightarrow 
 M=H(\RR)_{+}\cap \Gamma \backslash  X_{H}^+
$$
une application $\pi_{H,\alpha}: Z_{H}\longrightarrow M$.
Par ailleurs comme $M$ est un espace localement sym\'etrique 
hermitien, on dispose sur $M$ d'une m\'etrique k\"ahlerienne
et d'une mesure de probabilit\'e associ\'ee $\mu_{M}$. On a alors
\begin{equation}
(\pi_{H,\alpha})_{*}\mu_{H}= \mu_{M}.
\end{equation}
On remarque que $\pi_{H,\alpha}$ est la restriction \`a $Z_{H}$
de  l'application
\begin{equation}
\pi_{\alpha}: \Gamma \backslash G(\RR)^+\longrightarrow 
\Gamma \backslash X^+
\end{equation}
$$
\Gamma g\longrightarrow \Gamma g.\alpha
$$
Plus g\'en\'eralement pour tout $\beta\in X_{H}^+$, on a
$\pi_{\beta}(Z_{H})=M$ et $\pi_{\beta*}\mu_{H}=\mu_{M}$.

\medskip

\subsection{Sous-vari\'et\'es fortement sp\'eciales et retour vers 
les compacts.}

Nous gardons les notations de la partie pr\'ec\'edente.
La 
proposition suivante 
 est une forme faible d'un r\'esultat 
d\^u \`a Dani et Margulis
(\cite{DaMa} th\'eor\`eme 2)
\begin{prop}\label{prop4.5}
Il existe un ensemble compact $C$ de $\Gamma\backslash G(\RR)^+$
tel que pour tout sous-groupe unipotent \`a un param\`etre
$W=(u_{t})_{t\in\RR}$ et tout $g\in G(\RR)^+$, si
$$
\Gamma\backslash \Gamma g W\cap C=\emptyset
$$
alors il existe un $\QQ$-sous-groupe parabolique propre $P$ de $G_{\QQ}$
tel que 
$$gWg^{-1}\subset P(\RR).
$$
\end{prop}
On en d\'eduit tout d'abord:
\begin{lem}\label{lemme4.6}
Soient $F_{\QQ}$ un $\QQ$-sous-groupe semi-simple tel que
$F(\RR)^+\in \Hcal$ et $g\in G^+({\QQ})$ tels que 
$$
\Gamma \backslash \Gamma g F(\RR)^{+}\cap C=\emptyset.
$$
Il existe alors un $\QQ$-sous-groupe parabolique propre $P$ de $G_{\QQ}$
tel que $F_{\QQ}\subset P$.
\end{lem}
{\it Preuve.} Soit $L=L(F)$ le sous-groupe de $F(\RR)^{+}$ engendr\'e par 
les sous-groupes unipotents \`a un param\`etre de $G(\RR)^+$
qui sont contenus dans $F(\RR)^{+}$. On a alors $F_{\QQ}=MT(L)$ (voir 
la discussion apr\`es le lemme \ref{lem1}.
Fixons un sous-groupe \`a un param\`etre 
$$W=(u_{t})_{t\in \RR}\subset L$$
tel que $W$ n'est contenu dans aucun sous-groupe normal de $L(F)$.

Pour tout $h\in F(\RR)^+$
$$
\Gamma\backslash \Gamma ghW\cap C\subset 
\Gamma\backslash \Gamma gF(\RR)^+\cap C=\emptyset
$$
D'apr\`es la proposition \ref{prop4.5},
pour tout $h\in F(\RR)^+$, 
il existe un $\QQ$-parabolique propre $P_{h}$ tel que
$$
gh Wh^{-1}g^{-1}\subset P_{h}(\RR).
$$
Comme l'ensemble des paraboliques sur $\QQ$ est d\'enombrable,
il existe un $\QQ$-parabolique propre $P_{0}$ et un ensemble $A\subset 
F(\RR)^+$ de mesure
positive (pour la mesure de Haar sur $F(\RR)^{+}$) tel que pour 
tout $h\in A$:
$$
gh Wh^{-1}g^{-1}\subset P_{0}(\RR).
$$
Comme l'ensemble des $h\in F(\RR)^+$ tels que 
$gh Wh^{-1}g^{-1}\subset P_{0}(\RR)$ est un 
sous-ensemble de Zariski de $F(\RR)^+$
de mesure positive, par connexit\'e de $F(\RR)^+$ on en
d\'eduit que pour tout $h\in F(\RR)^+$
$$
ghWh^{-1}g^{-1}\subset P_{0\RR}.
$$
Comme $W$ n'est contenu dans aucun sous-groupe normal de $L$
on en d\'eduit que
$$
g L g^{-1}\subset P_{0\RR}.
$$
Soit $P=g^{-1}P_{0}g$, on a $L\subset P_{\RR}$ donc 
$$
F_{\QQ}=MT(L)\subset P.
$$
Ceci termine la preuve du lemme car $P$ est un $\QQ$-parabolique
propre.

\begin{lem}\label{Lemme4.5}
Il existe un compact $C'$ de $\Gamma \backslash X^+$ tel que si
$\Delta_{F} \backslash X_{F}^+$ est une sous-vari\'et\'e fortement 
sp\'eciale de type Shimura alors
\begin{equation}
\Delta_{F} \backslash X_{F}^+\cap C'\neq \emptyset
\end{equation}	 
\end{lem}	
{\it Preuve.} Soit $\Omega$ un compact de $G(\RR)^+$ contenant 
l'origine $e$ comme point interieur. On note
$$
C_{1}=C.\Omega=\{ c\omega, c\in C, \omega\in \Omega   \},
$$
c'est encore un compact de $\Gamma \backslash G(\RR)^{+}$.
Pour tout $x\in X^+$, on note comme pr\'ec\'edemment
$\pi_{x}$ l'application associ\'ee de $\Gamma \backslash G(\RR)^{+}$
vers $\Gamma \backslash X^+$. On fixe $x_{0}\in X^+$
et on note $C'=\pi_{x_{0}} (C_{1})$. On remarque que pour
tout $\omega\in \Omega$, si on note $x_{\omega}=\omega.x_{0}$
alors
$$
\pi_{x_{\omega}}(C)\subset C'.
$$

Soit $x\in X_{F}^+$ de sorte que $X_{F}^+=F(\RR)_{+}. x$.
Comme $G(\QQ)^+$ est dense dans $G(\RR)^+$, il existe
$\omega\in \Omega$ et $\gamma\in G(\QQ)^+$ tels que
$x=\gamma \omega.x_{0}$. Soit 
$F_{\gamma \QQ}=\gamma^{-1}F_{\QQ} \gamma$.
On a alors $X_{F}^+= F(\RR)_{+}\gamma.x_{\omega}$ et 
$$
\Delta_{F} \backslash X_{F}^+= \pi_{x_{\omega}} (\Gamma\backslash
\Gamma \gamma F_{\gamma}(\RR)_{+}).
$$
Si $\Delta_{F} \backslash X_{F}^+\cap C'=\emptyset$ alors \`a fortiori
$\Delta_{F} \backslash X_{F}^+\cap \pi_{x_{\omega}}(C)=\emptyset$
et finalement
$$
\Gamma\backslash
\Gamma \gamma F_{\gamma}(\RR)_{+}\cap C=\emptyset.
$$
Par le lemme \ref{lemme4.6}, on en d\'eduit que 
$F_{\gamma \QQ}\subset P$ pour un $\QQ$-sous-groupe parabolique propre $P$.
Ceci est impossible par la propri\'et\'e (b) des sous-vari\'et\'es 
fortement sp\'eciales.

\subsection{Preuve du th\'eor\`eme}

On dispose de tous les outils pour d\'emontrer le r\'esultat principal 
de ce texte:

\begin{teo} \label{teo4.3 }
	 Soient $(G,X)$ une donn\'ee de Shimura, $K\subset G(\AAA_{f})$ un
sous-groupe compact ouvert et $S=\Gamma \backslash X^{+}$
une composante irr\'eductible de $Sh_{K}(G,X)(\CC)$ pour un
$\Gamma=\Gamma_{\beta}$,  $\beta\in R_{G,K}$.  
Soit $S_{n}$ une suite de sous-vari\'et\'es fortement sp\'eciales de 
$S$.
Soit $\mu_{n}$ la mesure de probabilit\'e associ\'ee \`a $S_{n}$.
Il existe une sous-vari\'et\'e fortement sp\'eciale $M$ et une
sous-suite $\mu_{n_{k}}$ qui converge faiblement vers 
la mesure  $\mu_{M}$ canoniquement associ\'ee \`a $M$.
De plus $M$ contient $S_{n_{k}}$ pour tout $k$ assez grand.
\end{teo}

\begin{rem}
{\rm Sans la condition (a) des sous-vari\'et\'es fortement sp\'eciales
le th\'eor\`eme peut \^etre mis en d\'efaut, par exemple pour
des suites de tores $H_{n}$ d\'efinissant des points CM.
De m\^eme si (b) n'est pas v\'erifi\'e pour un groupe $H_{\QQ}$, on peut 
d'apr\`es la condition \'equivalente (b")   trouver
une suite $z_{n}\in Z_{G}(H)$ telle que l'image de $z_{n}$
dans $\Gamma \cap Z_{G}(H)\backslash Z_{G}(H)(\RR)$ n'a pas
de sous-suite convergente. Soit $\alpha: \SSS\rightarrow G_{\RR}$
se factorisant par $H_{\RR}$ et  $X_{n}$ la $H(\RR)$-classe de
conjuguaison de $z_{n}.\alpha$; alors $(H, X_{n})$ est une 
sous-donn\'ee de Shimura.  On peut v\'erifier que la suite de mesures
 canoniques $\mu_{n}$ 
sur la sous-vari\'et\'e sp\'eciale $\Gamma\cap 
H\backslash X_{n}^+$ n'a pas de sous-suite convergente.е}
\end{rem}

{\it Preuve}. On peut tout d'abord supposer que $G$ est adjoint.
Cela r\'esulte des d\'efinitions de sous-vari\'et\'es fortement
sp\'eciales en termes de la donn\'ee de Shimura adjointe
$(G^{ad},X^{ad})$ et de compatibilit\'es \'evidentes pour les mesures
canoniques des sous-vari\'etes fortement sp\'eciales de 
$S$ et de $$S^{ad}=Sh_{K^{ad}}(G^{ad},X^{ad})(\CC).$$ 

On peut supposer que les $S_{n}$ sont des sous-vari\'et\'es de type 
Shimura. En effet par le lemme \ref{lemme3.7}, en extrayant au besoin
une sous-suite, on peut supposer qu'il existe $\lambda\in R_{G,K}$
tel que $S_{n}$ est une composante irr\'eductible de $T_{\lambda}.S'_{n}$
pour une sous-vari\'et\'e fortement sp\'eciale $S'_{n}$ de type
Shimura. Le r\'esultat pour $S_{n}$ se d\'eduit alors de celui pour
$S'_{n}$.

Soit donc $S_{n}$ une suite de sous-vari\'et\'e fortement sp\'eciales
de type Shimura
de $S$.  
Soit $H_{n,\QQ}$ des $\QQ$-sous groupes associ\'es.  Soit 
$\alpha_{n}\in X^+$ tel que $MT(\alpha_{n})=T_{n}$
soit un $\QQ$-tore contenu dans $H_{n,\QQ}$е.
 On note alors $X_{n}^+$
la $H_{n}(\RR)_{+}$ classe de conjugaison de $\alpha_{n}$ et
$\Delta_{n}=\Gamma\cap H_{n}(\RR)_{+}$ de sorte que
$S_{n}= \Delta_{n}\backslash X_{n}^+$.

\begin{lem}\label{lemme4.3} 
	 Il existe une suite $\lambda_{n}\in \Gamma$ telle que si l'on note
$H_{n,\lambda_{n}е}е$ et	 $X^+_{n,\lambda_{n}е}е$ les conjugu\'es de
$H_{n}$ et $X^+_{n}е$ d\'efinis
par le proc\'ed\'e d\'ecrit dans la remarque \ref{remarque}, alors
en passant au besoin \`a une sous-suite,
Il existe une suite $\beta_{n}\in X_{n,\lambda_{n}е}^+$ qui converge
vers $\beta\in X^+$.
\end{lem}	 
D'apr\`es le lemme \ref{Lemme4.5}, il existe un compact $C'$ de 
$\Gamma\backslash X^+$ tel que
pour tout $n\in \NN$, 
$$
\Delta_{n}\backslash X_{n}^+\cap C'\neq \emptyset 
$$
Soit donc $t_{n}\in \Gamma_{n}\backslash X_{n}^+\cap C'$. On peut supposer
en passant au besoin \`a une sous-suite que $t_{n}\longrightarrow t\in 
C'$. Soit 
$$
\theta : X^+ \longrightarrow \Gamma \backslash X^+
$$
la projection. 
Les composantes irr\'eductibles des images inverses par $\theta$
de $S_{n}$ sont de la forme $\lambda.X_{H_{n}}^+=X^+_{H_{n,\lambda}}$ pour
un $\lambda\in \Gamma$. On peut donc en conjuguant au besoin $H_{n}$
par un $\lambda_{n}\in \Gamma$ choisir des 
 relev\'es 
par $\theta$ convenables $\beta_{n}\in X_{H_{n,\lambda_{n}е}}$
 de $t_{n}$ dans un domaine fondamental fixe pour l'action de $\Gamma$
 sur $X^+$. Alors la suite $\beta_{n}\longrightarrow \beta$ pour
 un relev\'e $\beta$ convenable de $t$.
 
 \medskip

On peut sans perte de g\'en\'eralit\'e supposer que 
$H_n=H_{n,\lambda_{n}е}е$.
On a vu 
que $H_{n,\QQ}(\RR)^+\in \Hcal$.
Pour 
tout $n$ on note $\mu'_{n}$ la mesure de $Q(\Omega,e)$
de support 
$$
\Gamma\backslash \Gamma H_{n}(\RR)^+.
$$
D'apr\`es la proposition \ref{prop2.5}, au besoin en passant \`a une
sous-suite, on peut supposer que $\mu'_{n}$ converge faiblement
vers une mesure $\mu'\in Q(\Omega,e)$. De plus pour tout $n$ assez
grand on a 
$$
supp(\mu_{n})= \Gamma\backslash \Gamma H_{n}(\RR)^+\subset supp(\mu').
$$
D'apr\`es la description  de $Q(\Omega,e)$ donn\'ee avant la 
proposition \ref{prop2.5} et le lemme \ref{lem1}, il existe 
$F\in \Hcal$ et un $\QQ$-sous-groupe alg\'ebrique $H_{\QQ}$
tel que $F=H_{\QQ}(\RR)^+$, $H_{\QQ}=MT(F)$ et $\mu'$
est la mesure $H(\RR)^+$-invariante de support
$supp(\mu)=\Gamma\backslash \Gamma H(\RR)^+$.
On en d\'eduit donc que 
$$
\Gamma\backslash \Gamma H_{n}(\RR)^+
\subset \Gamma\backslash \Gamma H(\RR)^+.
$$
On en d\'eduit que $Lie(H_{n}(\RR)^+)\subset Lie(H(\RR)^+)$ puis
par connexit\'e que 
$$H_{n}(\RR)^+\subset H(\RR)^+.$$
Finalement on obtient
$$
H_{n,\QQ}=MT(H_{n}(\RR)^+)\subset  H_{\QQ}=MT(H(\RR)^+).
$$
Pour tout $n$ assez grand, on a donc 
$$
T_{n}=MT(\alpha_{n})\subset H_{n,\QQ}\subset H_{\QQ},
$$
et la $H(\RR)_{+}$-classe de conjugaison  de $\alpha_{n}$
est ind\'ependante de $n$. On la note $X_{H}^+$.  Soit 
$\Delta_{H}=\Gamma \cap H(\RR)_{+}$.
D'apr\`es le lemme
\ref{lemme4.1}, $M=\Delta_{H}\backslash X_{H}^+$ est une 
sous-vari\'et\'e fortement sp\'eciale et pour tout $n$ assez grand
$$
S_{n}
\subset \Delta_{H}\backslash X_{H}^+.
$$
On termine la d\'emonstration 
de la mani\`ere suivante: pour tout $n\in \NN$, on a vu que
$\pi_{\beta_{n}*}\mu'_{n}=\mu_{n}$.  Comme $\beta_{n}\rightarrow 
\beta$, $\pi_{\beta_{n}}$ converge simplement et uniform\'ement sur
tout compacts vers $\pi_{\beta}$ et $\pi_{\beta*}\mu'=\mu_{M}$.
Soit  $f$ une fonction continue \`a support compact
sur $\Gamma\backslash X^+$.  On a 
$$
\mu_{n}(f)-\mu_{M}(f)= \mu'_{n}(f\pi_{\beta_{n}})- \mu'(f\pi_{\beta})
=\mu'_{n}(f\pi_{\beta_{n}})-\mu'_{n}(f\pi_{\beta})
+\mu'_{n}(f\pi_{\beta})- \mu'(f\pi_{\beta}).
$$
Comme $\mu'_{n}$ converge faiblement vers $\mu'$,
$\mu'_{n}(f\pi_{\beta})- \mu'(f\pi_{\beta})$ tend vers $0$.
Par la convergence uniforme sur les compacts de $\pi_{\beta_{n}}$ vers
$\pi_{\beta}$ et le fait que les $\mu_{n}$ sont des mesures de 
probabilit\'es, $\mu'_{n}(f\pi_{\beta_{n}})-\mu'_{n}(f\pi_{\beta})$
converge aussi vers $0$. On en d\'eduit donc que
$\mu_{n}(f)-\mu_{M}(f)$ converge vers $0$, donc que 
$\mu_{n}$ converge faiblement vers $\mu_{M}$.

\end{document}